\documentclass[10pt]{article} \usepackage{amssymb,amscd}
\usepackage[all]{xy} \usepackage[cp1255]{inputenc}
\usepackage[hebrew,english]{babel} \usepackage{amsmath}
\usepackage{xspace}
\usepackage{pdflscape}
\usepackage{array}
\makeatletter
\newcommand{\thickhline}{%
    \noalign {\ifnum 0=`}\fi \hrule height 1pt
    \futurelet \reserved@a \@xhline
}
\newcolumntype{"}{@{\hskip\tabcolsep\vrule width 1pt\hskip\tabcolsep}}
\makeatother
\makeindex
\author{Oz Ben-Shimol
 \\ \\
\small{Department of Mathematics,
Bar-Ilan University,} \\
\small{Ramat-Gan, Israel}}
\date{\small{\today}}
\title{Centers associated with the Borel subalgebra\\ of certain simple Lie algebras\thanks{Work is supported by ISF grant \#170/12.}}
\begin{document}
\maketitle
\newcommand{\gr}{\operatorname{gr}}
\newcommand{\si}{\operatorname{si}}
\newcommand{\ad}{\operatorname{ad}}
\newcommand{\tr}{\operatorname{trace}}
\newcommand{\depth}{\operatorname{depth}}
\newcommand{\Mod}{\operatorname{mod}}
\newcommand{\Ht}{\operatorname{ht}}
\newcommand{\Cl}{\operatorname{Cl}}
\newcommand{\Spec}{\operatorname{Spec}}
\newcommand{\Sl}{\operatorname{sl}}
\newcommand{\Span}{\operatorname{Span}}
\newcommand{\Char}{\operatorname{char}}
\newcommand{\pf}{\operatorname{pf}}
\newcommand{\Id}{\operatorname{Id}}
\newcommand{\Der}{\operatorname{Der}}
\newcommand{\cpf}{\operatorname{pf_{c}}}
\newcommand{\codim}{\operatorname{codim}}
\newcommand{\rest}[1]{\lvert_{#1}}

\newcommand{\ZZ}{{\mathbb Z}} \newcommand{\RR}{{\mathbb R}}
\newcommand{\NN}{{\mathbb N}} \newcommand{\QQ}{{\mathbb Q}}
\newcommand{\fg}{{\mathfrak{g}}} \newcommand{\HH}{{\mathbb H}}
\newcommand{\fn}{{\mathfrak{n}}}
\newcommand{\fhn}{{\mathfrak{\hat{n}}}} 
\newcommand{\fr}{{\mathfrak{r}}}
\newcommand{\fm}{{\mathfrak{M}}}
\newcommand{\DD}{{\mathfrak{D}}}
\newcommand{\KK}{{\mathfrak{K}}}
\newcommand{\BB}{{\mathfrak{b}}}
\newcommand{\FF}{{\mathbb F}}
\newcommand{\cL}{{\mathcal{L}}}
\newcommand{\ol}{\overline}
\newcommand{\eee}{\hfill$\Box$}

\thispagestyle{empty}
\baselineskip=0.55cm
\noindent
\thispagestyle{empty}
\begin{center}
    \large \textbf{Abstract}
\end{center}

We continue the study in Ben-Shimol [1],[2] and consider a Borel subalgebra $\BB$ and its nil radical $\fn$ of the simple Lie algebras of types $G_2$, $F_4$, $C_n$ over arbitrary field. 
Let $\cL\in\{\fn, \BB\}$.
We establish here explicit realizations of the center $Z(\cL)$ and semi-center $Sz(\cL)$ of the enveloping algebra, the Poisson center $S(\cL)^{\cL}$ and Poisson semi-center $S(\cL)^{\cL}_{\si}$ of the symmetric algebra. We describe their structure as commutative rings and establish isomorphisms $Z(\cL)\cong~S(\cL)^{\cL}$,  $Sz(\cL)\cong S(\cL)^{\cL}_{\si}$.
\subsection*{1. \ The $\mathbf{G_2}$ case}

\ \ \ \ \ An explicit basis of a particular Lie algebra of type $G_2$ over an algebraically closed field of characteristic zero, realized as a set of $7\times 7$ matrices, can be found in [7, pp. 103-106].
From this basis we take a part which forms a basis of a Borel subalgebra: \ 
$h_1=e_{22}-e_{33}-e_{55}+e_{66}$, \ $h_2=e_{33}-e_{44}-e_{66}+e_{77}$, \ $x_1=\sqrt{2}(e_{12}-e_{51})-(e_{37}-e_{46})$, $x_2=\sqrt{2}(e_{13}-e_{61})-(e_{45}-e_{27})$, \
$x_3=\sqrt{2}(e_{17}-e_{41})-(e_{53}-e_{62})$, \ $x_4=e_{23}-e_{65}$, \ $x_5=e_{42}-e_{57}$, \ $x_6=e_{43}-e_{67}$.
A customary notation for the elements $x_1,\ldots, x_6$ is $g_1$, $g_2$, $g_{-3}$, $g_{1,-2}$, $g_{3,-1}$, $g_{3,-2}$, respectively.
The table of Lie brackets is presented below: 
$$
\begin{array}{ccccccccccccccccc}
 & \vline width 1pt & h_1 & \vline & h_2 & \vline width 1pt & x_1 & \vline & x_2 & \vline & x_3 & \vline  & x_4 & \vline & x_5 & \vline & x_6  \\ \thickhline
 h_1 & \vline width 1pt & 0 & \vline & 0 &\vline  width 1pt & -x_1 &\vline & x_2 &\vline & 0 &\vline & 2x_4 &\vline & -x_5 &\vline & x_6   \\ \hline
 h_2 & \vline width 1pt & 0 & \vline  & 0 & \vline width 1pt &  0 & \vline & -x_2 & \vline & -x_3 & \vline & -x_4 & \vline & -x_5 & \vline & -2x_6  \\ \thickhline
x_1  & \vline width 1pt & x_1 & \vline & 0 & \vline width 1pt & 0 & \vline & 2x_3 & \vline & 3x_5 & \vline & x_2 & \vline & 0 & \vline & 0 \\ \hline
x_2 & \vline width 1pt & -x_2 & \vline & x_2 & \vline width 1pt & -2x_3 & \vline & 0 & \vline & 3x_6 & \vline & 0 & \vline & 0 & \vline & 0 \\ \hline
x_3  & \vline width 1pt & 0 & \vline & x_3 & \vline width 1pt & -3x_5 & \vline & -3x_6 & \vline & 0 & \vline & 0 & \vline & 0 & \vline & 0 \\ \hline
x_4 & \vline width 1pt & -2x_4 & \vline & x_4 & \vline width 1pt & -x_2 & \vline & 0 & \vline & 0 & \vline & 0 & \vline & -x_6 & \vline & 0 \\ \hline
x_5 & \vline width 1pt & x_5 & \vline & x_5 & \vline width 1pt & 0 & \vline & 0 & \vline & 0 & \vline & x_6 & \vline & 0 & \vline & 0 \\ \hline
x_6 & \vline width 1pt & -x_6 & \vline & 2x_6 & \vline width 1pt & 0 & \vline & 0 & \vline & 0 & \vline & 0 & \vline & 0 & \vline & 0 \\ \thickhline
\end{array}
$$

In the sequel we shall consider the Lie algebra $\BB$ defined by the table above, as a Lie algebra over an arbitrary field $F$ of characteristic $p$ where $p\neq 2,3$ ($p=0$ is considerable).
$\{h_1,h_2\}$ generates the Cartan subalgebra and $\{x_1,\ldots,x_6\}$ generates the nil radical $\fn$.
\\ \indent
Consider the symmetric algebra $S(\fn)$ as the polynomial algebra in the variables $x_1,\ldots,x_6$ over~$F$. 
Denote $c_1=x_6$, $c_2=3x_1x_6-3x_2x_5+x_{3}^{2}$.
It is easy to verify that $c_1,c_2\in S(\fn)^{\fn}$, that is, $\ad{x_i}(c_j)=0$ for all $i=1,\ldots,6$, \ $j=1,2$. 
\\ \indent
Assume $p>3$.
We denote by $S_{p}(\fn)$ the polynomial subalgebra of $S(\fn)$ which generated over $F$ by $c_1$ and the $p$-th powers $x_{i}^{p}$: \
$S_{p}(\fn)=F[x_{1}^{p},x_{2}^{p},x_{3}^{p},x_{4}^{p},x_{5}^{p},c_{1}]$.
Clearly, $S_{p}(\fn)\subseteq S(\fn)^{\fn}$.
Denote by $Q(A)$ the quotient field of a domain $A$. 
Then the field extension $Q( S_{p}(\fn)[c_2])/Q( S_{p}(\fn))$ is of degree $p$ because $c_2$ does not belong to $S_{p}(\fn)$ while $c_{2}^{p}$ does.
Consider the elements $v_2=-\frac{1}{3}x_3$, \ $v_3=\frac{1}{3}x_2$, \ $v_4=x_5$, \ $v_5=-x_4$.
For each $i,j$ such that $2\leq i\leq j\leq 5$ we have $\ad{v_i}(x_j)=\delta_{i,j}c_1$.
Therefore, $x_2\not\in S(\fn)^{\fn}[x_3,x_4,x_5]$, \ $x_3\not\in S(\fn)^{\fn}[x_4,x_5]$, \ $x_4\not\in S(\fn)^{\fn}[x_5]$, and of course $x_5\not\in S(\fn)^{\fn}$.
Together with the fact $x_{i}^{p}\in S(\fn)^{\fn}$ for each $i$ we deduce that the field extension $Q(S(\fn)^{\fn}[x_2,x_3,x_4,x_5])/Q(S(\fn)^{\fn})$ is of degree $p^{4}$.
The field $Q(S(\fn))=F(x_1,\ldots,x_6)$ is an extension of degree $p^5$ of $Q(S_{p}(\fn))$. 
By degree consideration $Q(S(\fn)^{\fn})=Q(S_{p}(\fn)[c_2])$, that is, the domains $S(\fn)^{\fn}$ and  $S_{p}(\fn)[c_2]$ have the same quotient field. 
\\ \indent  
The goal is to prove equality of the rings: $S(\fn)^{\fn}=S_{p}(\fn)[c_2]$.
The ring extension $S(\fn)^{\fn}/S_{p}(\fn)[c_2]$ is integral. 
Therefore, it suffices to prove that $S_{p}(\fn)[c_2]$ is normal (integrally closed in its quotient field), equivalently, $S_{p}(\fn)[c_2]$ satisfies $(S_1)$ and $(R_2)$ (see [8, p.183]).
Consider the polynomial ring $R=S_{p}(\fn)[t_2]$ of the one variable $t_2$.
The polynomial $f=t_{2}^{p}-c_{2}^{p}$ is a prime element of $R$ thus, $S_p(\fn)[c_2]\cong R/Rf$.
From [5, 2.1.28] it follows that $S_p(\fn)[c_2]$ is a Cohen-Macaulay ring, hence satisfies $(S_1)$.
To prove $(R_1)$, we have to show that if $P$ is an element of the singular locus of $R$ such that $f\in P$, then $\Ht{P}>2$. 
So let $P$ be such a prime.  
$\partial{f}/\partial{x_{1}^{p}}=-3x_{6}^{p}$, \ $\partial{f}/\partial{x_{2}^{p}}=-3x_{5}^{p}$.
Thus $P$ contains the prime ideal $R(x_{5}^{p},x_{6}^{p},f)$ which is of height $3$.  
We deduce that $S(\fn)^{\fn}=S_{p}(\fn)[c_2]$. Finally, if $p=0$ then  $S(\fn)^{\fn}=F[c_1,c_2]$ by [2, section 3].
\\ \indent
Let $Z(\fn)$ be the center of the enveloping algebra $U(\fn)$.
We shall use the same notation $x_i$, $i=1,\ldots,6$ for the basis of $\fn$, consider it as a Lie subalgebra of $U(\fn)$.  
The elements in $Z(\fn)$ correspond to $c_1$,$c_2$ will respectively denoted by $z_1$,$z_2$. 
So $z_1=x_6$, $z_2=3x_1x_6-3x_2x_5+x_{3}^{2}$.
Suppose $p>3$.
The analogous polynomial ring to $S_{p}(\fn)$ in $Z(\fn)$ is $Z_{p}(\fn)=F[x_{1}^{p},x_{2}^{p},x_{3}^{p},x_{4}^{p},x_{5}^{p},z_{1}]$. 
$S(\fn)$ is isomorphic to the graded algebra of $U(\fn)$, and we set $S(\fn)=\gr{U(\fn)}$.
In particular $c_i=\gr{z_i}$, $x_{i}^{p}=\gr{x_{i}^{p}}=(\gr{x_{i}})^{p}$.
Therefore $S(\fn)^{\fn}=\gr(Z_{p}(\fn)[z_2])\subseteq\gr{Z(\fn)}$.
The inclusion $\gr{Z(\fn)}\subseteq S(\fn)^{\fn}$ is trivial. 
Since $S(\fn)^{\fn}=S_{p}(\fn)[c_2]$, from [3, p.180, Prop. 10(ii), section 2.9] we have $Z(\fn)=Z_{p}(\fn)[z_2]$.
\\ \indent
Let $\varphi:R\to Z(\fn)$ be the $F$-algebra epimorphism defined by $\varphi(x_{i}^{p})=x_{i}^{p}$, \ $\varphi(c_1)=z_1$, \
$\varphi(t_2)=z_{2}$. Obviously, $Rf\subseteq\ker\varphi$. Hence $Z(\fn)$ is a homomorphic image of $R/Rf$. 
The rings $Z(\fn)$ and $R/Rf$ are both domains with equal Krull dimension, hence $Z(\fn)\cong R/Rf$. 
We deduce that $Z(\fn)\cong S(\fn)^{\fn}$. 
The rings $Z(\fn), S(\fn)^{\fn}$ are hypersurfaces.
\\ \indent
If $p=0$, we have $S(\fn)^{\fn}=F[c_1,c_2]=F[\gr{z_1},\gr{z_2}]=\gr{F[z_1,z_2]}\subseteq\gr{Z(\fn)}$ and $\gr{Z(\fn)}\subseteq S(\fn)^{\fn}$.
Therefore $Z(\fn)=F[z_1,z_2]$ and $Z(\fn)\cong S(\fn)^{\fn}$ as polynomial algebras in two variables. 
We should remark here that an isomorphism $Z(\fn)\cong S(\fn)^{\fn}$ where $p=0$ in known [6, Proposition 4.8.12]. 
\\ \indent
Let summarize the main results that presented until now: 
\vskip 0.1cm\noindent 
\textbf{1.1. Theorem.} \textbf{a.} \ \textit{Suppose} $\mathit{p>3}$. \textit{Then} 
$$
	\mathit{S(\fn)^{\fn}=S_{p}(\fn)[c_2]\cong S_{p}(\fn)[t_2]/(t_{2}^{p}-c_{2}^{p})\cong Z(\fn)=Z_{p}(\fn)[z_2]}.
$$
\textit{In particular,} $\mathit{S(\fn)^{\fn}}$, $\mathit{Z(\fn)}$ \textit{are hypersurfaces.} \\ \indent
\textbf{b.} \ \textit{For} $\mathit{p=0}$, \ $\mathit{S(\fn)^{\fn}=F[c_1,c_2]\cong Z(\fn)=F[z_1,z_2]}$. \\
\textit{In particular,} $\mathit{S(\fn)^{\fn}}$, $\mathit{Z(\fn)}$ \textit{are polynomial rings in two variables.}
\vskip 0.2cm \indent
Let $S(\BB)_{\si}^{\BB}$ be the Poisson semi-center of $S(\BB)$. 
Assume $p>3$.
The linear transformations $\ad{x_i}, \ad{h_{j}}:S(\BB)\to S(\BB)$ satisfy the split equation $X^{p}-X=0$ over $F$. 
Together with the fact $\fn=[\BB,\BB]$, we have $S(\BB)_{\si}^{\BB}=S(\BB)^{\fn}$ (see [2, section 4]). 
Denote $S_{p}(\BB)=S_{p}(\fn)[h_{1}^{p},h_{2}^{p}]$. 
$S_{p}(\BB)$ is a polynomial subalgebra of $S(\BB)$.
Clearly, $S_{p}(\BB)[c_2]\subseteq S(\BB)^{\fn}$.
The field extension $Q(S_{p}(\BB)[c_2])/Q( S_{p}(\BB))$ is of degree $p$.
The field extension $Q(S(\BB)^{\fn}[x_2,x_3,x_4,x_5,h_1,h_2])/Q(S(\BB)^{\fn})$ is of degree $p^{6}$.
Indeed, as before, using the $v_{i}$'s we get that the field extension $Q(S(\BB)^{\fn}[x_2,x_3,x_4,x_5])/Q(S(\BB)^{\fn})$ is of degree $p^{4}$.
Now, $\ad{c_1}(h_1)=-c_1$ while $\ad{c_1}(x_{i})=0$, and $\ad{c_2}(h_2)=2c_2$ while $\ad{c_2}(x_{i})=\ad{c_2}(h_{1})=0$. 
Therefore, the field extension $Q(S(\BB)^{\fn}[x_2,x_3,x_4,x_5,h_1,h_2])/Q(S(\BB)^{\fn}[x_2,x_3,x_4,x_5])$ is of degree $p^{2}$.
The field $Q(S(\BB))=F(x_1,\ldots,x_6,h_1,h_2)$ is an extension of degree $p^7$ of $Q(S_{p}(\BB))$. 
By degree consideration $Q(S(\BB)^{\fn})=Q(S_{p}(\BB)[c_2])$, that is, the domains $S(\BB)^{\fn}$ and  $S_{p}(\BB)[c_2]$ have the same quotient field.  
\\ \indent
Identical arguments we applied to $S(\fn)^{\fn}$ yield $S(\BB)^{\fn}=S_{p}(\BB)[c_2]\cong S_{p}(\BB)[t_2]/(t_{2}^{p}-c_{2}^{p})$, where $t_2$ is transcendental over $S_{p}(\BB)$. 
Also, if $p=0$ then $S(\BB)^{\fn}=F[c_1,c_2]$.  
\\ \indent
Let $Sz(\BB)$ be the semi-center of $U(\BB)$. 
By [4, Proposition 2.1], $Sz(\BB)$ is commutative.
We shall use the same notation $h_1$, $h_2$ for the basis of the Cartan subalgebra of $\BB$, consider it as a Lie subalgebra of $U(\BB)$.
Assume $p>3$.
The analogous polynomial ring to $S_{p}(\BB)$ in $Sz(\BB)$ is $Z_{p}(\BB)=Z_{p}(\fn)[h_{1}^{p}-h_{1},h_{2}^{p}-h_{2}]$ ($x_{i}^{p}$ and $h_{j}^{p}-h_{j}$ are central (weight $0$) while $z_1$ is semi-central with non zero weight). 
$S(\BB)$ is isomorphic to the graded algebra of $U(\BB)$, and we set $S(\BB)=\gr{U(\BB)}$.
In particular $h_{i}^{p}=\gr(h_{i}^{p}-h_{i})$.
Therefore $S(\BB)_{\si}^{\BB}=S(\BB)^{\fn}=\gr(Z_{p}(\BB)[z_2])\subseteq\gr{Sz(\BB)}$ ($z_2$ is clearly semi-central with non zero weight).
The inclusion $\gr{Sz(\BB)}\subseteq S(\BB)_{\si}^{\BB}$ is trivial. 
Since $S(\BB)_{\si}^{\BB}=S_{p}(\BB)[c_2]$, from [3, p.180, Prop. 10(ii), section 2.9] we have $Sz(\BB)=Z_{p}(\BB)[z_2]$.
\\ \indent
Let $\varphi:S_{p}(\BB)[t_2]\to Sz(\BB)$ be the $F$-algebra epimorphism defined by $\varphi(x_{i}^{p})=x_{i}^{p}$, \ $\varphi(h_{j}^{p})=h_{j}^{p}-h_{j}$ \ $\varphi(c_1)=z_1$, \
$\varphi(t_2)=z_{2}$. 
Obviously, $t_{2}^{p}-c_{2}^{p}\in\ker\varphi$. 
Hence $Sz(\BB)$ is a homomorphic image of $S_{p}(\BB)[t_{2}]/(t_{2}^{p}-c_{2}^{p})$. 
The rings $Sz(\BB)$ and $S_{p}(\BB)[t_2]/(t_{2}^{p}-c_{2}^{p})$ are both domains with equal Krull dimension, hence $Sz(\BB)\cong S_{p}(\BB)[t_{2}]/(t_{2}^{p}-c_{2}^{p})$. 
We deduce that $Sz(\BB)\cong S(\BB)_{\si}^{\BB}$. 
The rings $Sz(\BB), S(\BB)_{\si}^{\BB}$ are hypersurfaces.
\\ \indent
If $p=0$, we have $S(\BB)_{\si}^{\BB}=F[c_1,c_2]=F[\gr{z_1},\gr{z_2}]=\gr{F[z_1,z_2]}\subseteq\gr{Sz(\BB)}$ and $\gr{Sz(\BB)}\subseteq S(\BB)_{\si}^{\BB}$.
Therefore $Sz(\BB)=F[z_1,z_2]$ and $Sz(\BB)\cong S(\BB)_{\si}^{\BB}$ as polynomial algebras in two variables. 
We should remark here that an isomorphism $Sz(\BB)\cong S(\BB)_{\si}^{\BB}$ where $p=0$ and $F$ is algebraically closed is known (see [9]). 
\\ \indent
Let summarize the results for the semi-centers: 
\vskip 0.1cm\noindent 
\textbf{1.2. Theorem.} \textbf{a.} \ \textit{Suppose} $\mathit{p>3}$. \textit{Then} 
$$
	\mathit{S(\BB)_{\si}^{\BB}=S_{p}(\BB)[c_2]\cong S_{p}(\BB)[t_{2}]/(t_{2}^{p}-c_{2}^{p})\cong Sz(\BB)=Z_{p}(\BB)[z_2]}.
$$
\textit{The rings} $\mathit{S(\BB)_{\si}^{\BB}}$, $\mathit{Sz(\BB)}$ \textit{are hypersurfaces.} \\ \indent
\textbf{b.} \ \textit{For} $\mathit{p=0}$, \ $\mathit{S(\BB)_{\si}^{\BB}=F[c_1,c_2]\cong Sz(\BB)=F[z_1,z_2]}$. \\
\textit{Therefore,} $\mathit{S(\BB)_{\si}^{\BB}=S(\fn)^{\fn}}$, \ $\mathit{Sz(\BB)=Z(\fn)}$
\textit{and the both rings} $\mathit{S(\BB)_{\si}^{\BB}}$, $\mathit{Sz(\BB)}$ \textit{are polynomial rings in two variables.}
\vskip 0.2cm \indent
equipped with the semi-centers we can find the Poisson center $S(\BB)^{\BB}$ of $S(\BB)$ and the center $Z(\BB)$ of $U(\BB)$ (reversing the approach in [2]).
\\ \indent
Suppose $p>3$.
Clearly, $Q(F[x_{1}^{p},\ldots,x_{6}^{p},h_{1}^{p},h_{2}^{p}])\subseteq Q(S(\BB)^{\BB})$.
Since $\ad{h_2}(c_2)=-c_{2}$, \ $\ad{h_1}(c_2)=0$, $\ad{h_1}(c_1)=c_1$, the extension $Q(S(\BB)^{\BB}[c_1,c_2])/Q(S(\BB)^{\BB})$ is of degree $p^{2}$. 
But $S(\BB)^{\BB}[c_1,c_2]\subseteq S(\BB)_{\si}^{\BB}=S_{p}(\BB)[c_2]$, and $Q(S_{p}(\BB)[c_2])$ is of degree $p^{2}$ over $Q(F[x_{1}^{p},\ldots,x_{6}^{p},h_{1}^{p},h_{2}^{p}])$.
By degree consideration we must have $Q(S(\BB)^{\BB})=Q(F[x_{1}^{p},\ldots,x_{6}^{p},h_{1}^{p},h_{2}^{p}])$, hence $S(\BB)^{\BB}=F[x_{1}^{p},\ldots,x_{6}^{p},h_{1}^{p},h_{2}^{p}]$. 
That is, $S(\BB)^{\BB}$ is generated over $F$ by the $p$-th powers of the basis elements of $\BB$.
For $p=0$ we of course have $S(\BB)^{\BB}=F$ by [2, section 3].
\\ \indent
Applying grading consideration one get $Z(\BB)=F[x_{1}^{p},\ldots,x_{6}^{p},h_{1}^{p}-h_1,h_{2}^{p}-h_2]$ for $p>3$ and $Z(\BB)=F$ for $p=0$. 
\\ \indent
So we have
\vskip 0.1cm\noindent 
\textbf{1.3. Theorem.} \textbf{a.} \ \textit{Suppose} $\mathit{p>3}$. \textit{Then} 
$$
\mathit{S(\BB)^{\BB}=F[x_{1}^{p},\ldots,x_{6}^{p},h_{1}^{p},h_{2}^{p}])\cong Z(\BB)=F[x_{1}^{p},\ldots,x_{6}^{p},h_{1}^{p}-h_1,h_{2}^{p}-h_2]},
$$
\textit{and are polynomial rings in eight variables.} \\ \indent
\textbf{b.} \ \textit{For} $\mathit{p=0}$, \ $\mathit{S(\BB)^{\BB}=Z(\BB)=F}$. 
\subsection*{2. \ The $\mathbf{F_4}$ case}

\ \ \ \ \ The positive roots and a corresponding basis of a nil radical $\fn$ of the Lie algebras of type $F_4$, over an algebraically closed field of characteristic zero, is presented below (the list of positive roots is taken from [10, page 274]):
$$
\begin{array}{llllll}
(1,0,0,0) & x_{1} & = & a & &  \\
(0,1,0,0) & x_{2} & = & b & &  \\
(0,0,1,0) & x_{3} & = & c & &  \\
(0,0,0,1) & x_{4} & = & d & &  \\
(1,1,0,0) & x_{5} & = & [x_{1},x_{2}] & = & [a,b] \\
(0,1,1,0) & x_{6} & = & [x_{2},x_{3}] & = & [b,c] \\
(0,0,1,1) & x_{7} & = & [x_{3},x_{4}] & = & [c,d] \\
(1,1,1,0) & x_{8} & = & [x_{1},x_{6}] & = & [a,[b,c]] \\
(0,1,2,0) & x_{9} & = & [x_{3},x_{6}] & = & [c,[b,c]] \\
(0,1,1,1) & x_{10} & = & [x_{2},x_{7}] & = & [b,[c,d]] \\
(1,1,2,0) & x_{11} & = & [x_{1},x_{9}] & = & [a,[c,[b,c]]] \\
(1,1,1,1) & x_{12} & = & [x_{1},x_{10}] & = & [a,[b,[c,d]]] \\
(0,1,2,1) & x_{13} & = & [x_{3},x_{10}] & = & [c,[b,[c,d]]] \\
(1,2,2,0) & x_{14} & = & [x_{2},x_{11}] & = & [b,[a,[c,[b,c]]]] \\
(1,1,2,1) & x_{15} & = & [x_{1},x_{13}] & = & [a,[c,[b,[c,d]]]] \\
(0,1,2,2) & x_{16} & = & [x_{4},x_{13}] & = & [d,[c,[b,[c,d]]]] \\
(1,2,2,1) & x_{17} & = & [x_{4},x_{14}] & = & [d,[b,[a,[c,[b,c]]]]] \\
(1,1,2,2) & x_{18} & = & [x_{1},x_{16}] & = & [a,[d,[c,[b,[c,d]]]]] \\
(1,2,3,1) & x_{19} & = & [x_{3},x_{17}] & = & [c,[d, [b,[a,[c,[b,c]]]]]] \\
(1,2,2,2) & x_{20} & = & [x_{4},x_{17}] & = & [d,[d,[b,[a,[c,[b,c]]]]]] \\
(1,2,3,2) & x_{21} & = & [x_{4},x_{19}] & = & [d,[c,[d, [b,[a,[c,[b,c]]]]]]] \\
(1,2,4,2) & x_{22} & = & [x_{3},x_{21}] & = & [c,[d,[c,[d, [b,[a,[c,[b,c]]]]]]]] \\
(1,3,4,2) & x_{23} & = & [x_{2},x_{22}] & = & [b,[c,[d,[c,[d, [b,[a,[c,[b,c]]]]]]]]] \\
(2,3,4,2) & x_{24} & = & [x_{1},x_{23}] & = & [a,[b,[c,[d,[c,[d, [b,[a,[c,[b,c]]]]]]]]]]
\end{array}
$$
\vskip 0.2cm\noindent
Cartan matrix:
$$
\left(
\begin{array}{rrrr}
	2 & -1 & 0 & 0 \\
	-1 & 2 & -2 & 0 \\
	0 & -1 & 2 & -1 \\
	0 & 0 & -1 & 2
\end{array}	
\right)
$$
\begin{landscape}
\ \ \ \ \ \ \ A fixed basis of the Cartan subalgebra: \ $h_1,h_2,h_3,h_4$. \vskip 0.1cm
\ \ \ \ \ \ \ The  Borel subalgebra generated by the $h_{i}$'s and the $x_{j}$'s will denoted by $\BB$.
\vskip 0.2cm\noindent
\ \ \ \ \ \ \ \ \ \ \ The table of Lie brackets is presented below:
$$
\begin{array}{ccccccccccccccccccccccccc}
 & \vline width 1pt & x_1 & \vline & x_2 & \vline & x_3 & \vline & x_4 & \vline & x_5 & \vline  & x_6 & \vline & x_7 & \vline & x_8  & \vline & x_9 & \vline & x_{10} & \vline & x_{11} & \vline & x_{12} \\ \thickhline
 h_1 & \vline width 1pt & 2x_1 & \vline & -x_2 &\vline & 0 &\vline & 0 &\vline & x_5 &\vline & -x_6 &\vline & 0 &\vline & x_8 & \vline & -x_9 & \vline & -x_{10} & \vline & x_{11} & \vline & x_{12}  \\ \hline
h_2 & \vline width 1pt & -x_1 & \vline  & 2x_2 & \vline & -x_3 & \vline & 0 & \vline & x_5 & \vline & x_6 & \vline & -x_7 & \vline & 0 & \vline & 0 & \vline & x_{10} & \vline & -x_{11} & \vline & 0  \\ \hline
 h_3 & \vline width 1pt & 0     & \vline  & -2x_2 & \vline & 2x_3 & \vline & -x_4 & \vline & -2x_5 & \vline & 0 & \vline & x_7 & \vline & 0 & \vline & 2x_9 & \vline & -x_{10} & \vline & 2x_{11} & \vline & -x_{12}  \\ \hline
 h_4 & \vline width 1pt & 0     & \vline  & 0 & \vline & -x_3 & \vline & 2x_4 & \vline & 0 & \vline & -x_6 & \vline & x_7 & \vline & -x_8 & \vline & -2x_9 & \vline & x_{10} & \vline & -2x_{11} & \vline & x_{12}  \\ \thickhline
 x_1 & \vline width 1pt & 0     & \vline  & x_5 & \vline & 0 & \vline & 0 & \vline & 0 & \vline & x_8 & \vline & 0 & \vline & 0 & \vline & x_{11} & \vline & x_{12} & \vline & 0 & \vline & 0  \\ \hline
 x_2 & \vline width 1pt & -x_5     & \vline  & 0 & \vline & x_6 & \vline & 0 & \vline & 0 & \vline & 0 & \vline & x_{10} & \vline & 0 & \vline & 0 & \vline & 0 & \vline & x_{14} & \vline & 0  \\ \hline
 x_3 & \vline width 1pt & 0     & \vline  & -x_6 & \vline & 0 & \vline & x_7 & \vline & -x_8 & \vline & x_9 & \vline & 0 & \vline & x_{11} & \vline & 0 & \vline & x_{13} & \vline & 0 & \vline & x_{15}  \\ \hline
 x_4 & \vline width 1pt & 0     & \vline  & 0 & \vline & -x_7 & \vline & 0 & \vline & 0 & \vline & -x_{10} & \vline & 0 & \vline & -x_{12} & \vline & -2x_{13} & \vline & 0 & \vline & -2x_{15} & \vline & 0  \\ \hline
 x_5 & \vline width 1pt & 0     & \vline  & 0 & \vline & x_8 & \vline & 0 & \vline & 0 & \vline & 0 & \vline & x_{12} & \vline & 0 & \vline & -x_{14} & \vline & 0 & \vline & 0 & \vline & 0  \\ \hline
 x_6 & \vline width 1pt & -x_8     & \vline  & 0 & \vline & -x_9 & \vline & x_{10} & \vline & 0 & \vline & 0 & \vline & -x_{13} & \vline & x_{14} & \vline & 0 & \vline & 0 & \vline & 0 & \vline & -\frac{1}{2}x_{17}  \\ \hline
 x_7 & \vline width 1pt & 0     & \vline  & -x_{10} & \vline & 0 & \vline & 0 & \vline & -x_{12} & \vline & x_{13} & \vline & 0 & \vline & x_{15} & \vline & 0 & \vline & -x_{16} & \vline & 0 & \vline & -x_{18}  \\ \hline
 x_8 & \vline width 1pt & 0     & \vline  & 0 & \vline & -x_{11} & \vline & x_{12} & \vline & 0 & \vline & -x_{14} & \vline & -x_{15} & \vline & 0 & \vline & 0 & \vline & \frac{1}{2}x_{17} & \vline & 0 & \vline & 0  \\ \hline
 x_9 & \vline width 1pt & -x_{11}     & \vline  & 0 & \vline & 0 & \vline & 2x_{13} & \vline & x_{14} & \vline & 0 & \vline & 0 & \vline & 0 & \vline & 0 & \vline & 0 & \vline & 0 & \vline & -x_{19}  \\ \hline
 x_{10} & \vline width 1pt & -x_{12}     & \vline  & 0 & \vline & -x_{13} & \vline & 0 & \vline & 0 & \vline & 0 & \vline & x_{16} & \vline & -\frac{1}{2}x_{17} & \vline & 0 & \vline & 0 & \vline & -x_{19} & \vline & \frac{1}{2}x_{20}  \\ \hline
 x_{11} & \vline width 1pt & 0     & \vline  & -x_{14} & \vline & 0 & \vline & 2x_{15} & \vline & 0 & \vline & 0 & \vline & 0 & \vline & 0 & \vline & 0 & \vline & x_{19} & \vline & 0 & \vline & 0  \\ \hline
 x_{12} & \vline width 1pt & 0     & \vline  & 0 & \vline & -x_{15} & \vline & 0 & \vline & 0 & \vline & \frac{1}{2}x_{17} & \vline & x_{18} & \vline & 0 & \vline & x_{19} & \vline & -\frac{1}{2}x_{20} & \vline & 0 & \vline & 0  \\ \hline
 x_{13} & \vline width 1pt & -x_{15}     & \vline  & 0 & \vline & 0 & \vline & -x_{16} & \vline & -\frac{1}{2}x_{17} & \vline & 0 & \vline & 0 & \vline & \frac{1}{2}x_{19} & \vline & 0 & \vline & 0 & \vline & 0 & \vline & \frac{1}{2}x_{21}  \\ \hline
 x_{14} & \vline width 1pt & 0     & \vline  & 0 & \vline & 0 & \vline & -x_{17} & \vline & 0 & \vline & 0 & \vline & -x_{19} & \vline & 0 & \vline & 0 & \vline & 0 & \vline & 0 & \vline & 0  \\ \hline
 x_{15} & \vline width 1pt & 0     & \vline  & \frac{1}{2}x_{17} & \vline & 0 & \vline & -x_{18} & \vline & 0 & \vline & -\frac{1}{2}x_{19} & \vline & 0 & \vline & 0 & \vline & 0 & \vline & -\frac{1}{2}x_{21} & \vline & 0 & \vline & 0  \\ \hline
 x_{16} & \vline width 1pt & -x_{18}     & \vline  & 0 & \vline & 0 & \vline & 0 & \vline & -\frac{1}{2}x_{20} & \vline & 0 & \vline & 0 & \vline & x_{21} & \vline & 0 & \vline & 0 & \vline & x_{22} & \vline & 0  \\ \hline
 x_{17} & \vline width 1pt & 0     & \vline  & 0 & \vline & -x_{19} & \vline & -x_{20} & \vline & 0 & \vline & 0& \vline & -x_{21} & \vline & 0 & \vline & 0 & \vline & 0 & \vline & 0 & \vline & 0  \\ \hline
 x_{18} & \vline width 1pt & 0     & \vline  & \frac{1}{2}x_{20} & \vline & 0 & \vline & 0 & \vline & 0 & \vline & -x_{21} & \vline & 0 & \vline & 0 & \vline & -x_{22} & \vline & 0 & \vline & 0 & \vline &0  \\ \hline
 x_{19} & \vline width 1pt & 0     & \vline  & 0 & \vline & 0 & \vline & -x_{21} & \vline & 0 & \vline & 0 & \vline & -x_{22} & \vline & 0 & \vline & 0 & \vline & -x_{23} & \vline & 0 & \vline & -x_{24}  \\ \hline
 x_{20} & \vline width 1pt & 0     & \vline  & 0 & \vline & -2x_{21} & \vline & 0 & \vline & 0 & \vline & 0 & \vline & 0 & \vline & 0 & \vline & 2x_{23} & \vline & 0 & \vline & 2x_{24} & \vline & 0  \\ \hline
 x_{21} & \vline width 1pt & 0     & \vline  & 0 & \vline & -x_{22} & \vline & 0 & \vline & 0 & \vline & -x_{23} & \vline & 0 & \vline & -x_{24} & \vline & 0 & \vline & 0 & \vline & 0 & \vline & 0  \\ \hline
 x_{22} & \vline width 1pt & 0     & \vline  & -x_{23} & \vline & 0 & \vline & 0 & \vline & -x_{24} & \vline & 0 & \vline & 0 & \vline & 0 & \vline & 0 & \vline & 0 & \vline & 0 & \vline & 0  \\ \hline
 x_{23} & \vline width 1pt & -x_{24}     & \vline  & 0 & \vline & 0 & \vline & 0& \vline & 0 & \vline & 0 & \vline & 0 & \vline & 0 & \vline & 0 & \vline & 0 & \vline & 0 & \vline & 0  \\ \hline
 x_{24} & \vline width 1pt & 0     & \vline  & 0 & \vline & 0 & \vline & 0 & \vline & 0 & \vline & 0 & \vline & 0 & \vline & 0 & \vline & 0 & \vline & 0 & \vline & 0 & \vline & 0  \\ \thickhline
\end{array}
$$
\end{landscape}
\begin{landscape}
$$
\begin{array}{cccccccccccccccccccccccccc}
 & \vline width 1pt & x_{13} & \vline & x_{14} & \vline & x_{15} & \vline & x_{16} & \vline & x_{17} & \vline  & x_{18} & \vline & x_{19} & \vline & x_{20}  & \vline & x_{21} & \vline & x_{22} & \vline & x_{23} & \vline & x_{24} & \vline width 1pt \\ \thickhline
 h_1 & \vline width 1pt & -x_{13} & \vline & 0 &\vline & x_{15} &\vline & -x_{16} &\vline & 0 &\vline & x_{18} &\vline & 0 &\vline & 0 & \vline & 0 & \vline & 0 & \vline & -x_{23} & \vline & x_{24} & \vline width 1pt  \\ \hline
 h_2 & \vline width 1pt & 0 & \vline  & x_{14} & \vline & -x_{15} & \vline & 0 & \vline & x_{17} & \vline & -x_{18} & \vline & 0 & \vline & x_{20} & \vline & 0 & \vline & -x_{22} & \vline & x_{23} & \vline & 0 & \vline width 1pt \\ \hline
h_3 & \vline width 1pt & x_{13} & \vline & 0 & \vline & x_{15} & \vline & 0 & \vline & -x_{17} & \vline & 0 & \vline & x_{19} & \vline & -2x_{20} & \vline & 0 & \vline & 2x_{22} & \vline & 0 & \vline & 0 & \vline width 1pt \\ \hline
h_4 & \vline width 1pt & 0     & \vline  & -2x_{14}& \vline & 0 & \vline & 2x_{16} & \vline & 0 & \vline & 2x_{18} & \vline & -x_{19} & \vline & 2x_{20} & \vline & x_{21} & \vline & 0 & \vline & 0 & \vline & 0  & \vline width 1pt\\ \thickhline
 x_1 & \vline width 1pt & x_{15}     & \vline  & 0 & \vline & 0 & \vline & x_{18} & \vline & 0 & \vline & 0 & \vline & 0 & \vline & 0 & \vline & 0 & \vline & 0 & \vline & x_{24} & \vline & 0 & \vline width 1pt \\ \hline
 x_2 & \vline width 1pt & 0  & \vline  & 0 & \vline &  -\frac{1}{2}x_{17} & \vline & 0 & \vline & 0 & \vline &  -\frac{1}{2}x_{20} & \vline & 0 & \vline & 0 & \vline & 0 & \vline & x_{23} & \vline & 0 & \vline & 0 & \vline width 1pt \\ \hline
 x_3 & \vline width 1pt & 0     & \vline  & 0 & \vline & 0 & \vline & 0 & \vline & x_{19} & \vline & 0 & \vline & 0 & \vline & 2x_{21} & \vline & x_{22} & \vline & 0 & \vline & 0 & \vline & 0 & \vline width 1pt \\ \hline
 x_4 & \vline width 1pt & x_{16}  & \vline  & x_{17} & \vline & x_{18} & \vline & 0 & \vline & x_{20} & \vline & 0 & \vline & x_{21} & \vline & 0 & \vline & 0 & \vline & 0 & \vline & 0 & \vline & 0 & \vline width 1pt\\ \hline
 x_5 & \vline width 1pt &  \frac{1}{2}x_{17} & \vline  & 0 & \vline & 0 & \vline &  \frac{1}{2}x_{20} & \vline & 0 & \vline & 0 & \vline & 0 & \vline & 0 & \vline & 0 & \vline & x_{24} & \vline & 0 & \vline & 0 & \vline width 1pt \\ \hline
 x_6 & \vline width 1pt & 0 & \vline  & 0 & \vline &  \frac{1}{2}x_{19} & \vline & 0 & \vline & 0 & \vline & x_{21} & \vline & 0 & \vline & 0 & \vline & x_{23} & \vline & 0 & \vline & 0 & \vline & 0 & \vline width 1pt \\ \hline
 x_7 & \vline width 1pt & 0     & \vline  & x_{19} & \vline & 0 & \vline & 0 & \vline & x_{21} & \vline & 0 & \vline & x_{22} & \vline & 0 & \vline & 0 & \vline & 0 & \vline & 0 & \vline & 0 & \vline width 1pt \\ \hline
 x_8 & \vline width 1pt &  -\frac{1}{2}x_{19} & \vline  & 0 & \vline & 0 & \vline & -x_{21} & \vline & 0 & \vline & 0 & \vline & 0 & \vline & 0 & \vline & x_{24} & \vline & 0 & \vline & 0 & \vline & 0  & \vline width 1pt\\ \hline
 x_9 & \vline width 1pt & 0  & \vline  & 0 & \vline & 0 & \vline & 0 & \vline & 0 & \vline & x_{22} & \vline & 0 & \vline & -2x_{23} & \vline & 0 & \vline & 0 & \vline & 0 & \vline & 0 & \vline width 1pt \\ \hline
 x_{10} & \vline width 1pt & 0 & \vline  & 0 & \vline &  \frac{1}{2}x_{21} & \vline & 0 & \vline & 0 & \vline & 0 & \vline & x_{23} & \vline & 0 & \vline & 0 & \vline & 0 & \vline & 0 & \vline & 0 & \vline width 1pt \\ \hline
 x_{11} & \vline width 1pt & 0     & \vline  & 0 & \vline & 0 & \vline & -x_{22} & \vline & 0 & \vline & 0 & \vline & 0 & \vline & -2x_{24} & \vline & 0 & \vline & 0 & \vline & 0 & \vline & 0 & \vline width 1pt \\ \hline
 x_{12} & \vline width 1pt &  -\frac{1}{2}x_{21} & \vline  & 0 & \vline & 0 & \vline & 0 & \vline & 0 & \vline & 0 & \vline & x_{24} & \vline & 0 & \vline & 0 & \vline & 0 & \vline & 0 & \vline & 0 & \vline width 1pt \\ \hline
 x_{13} & \vline width 1pt & 0 & \vline  & 0 & \vline & \frac{1}{2}x_{22} & \vline & 0 & \vline & -x_{23} & \vline & 0 & \vline & 0 & \vline & 0 & \vline & 0 & \vline & 0 & \vline & 0 & \vline & 0 & \vline width 1pt \\ \hline
 x_{14} & \vline width 1pt & 0     & \vline  & 0 & \vline & 0 & \vline & -x_{23} & \vline & 0 & \vline & -x_{24} & \vline & 0 & \vline & 0 & \vline & 0 & \vline & 0 & \vline & 0 & \vline & 0 & \vline width 1pt \\ \hline
 x_{15} & \vline width 1pt &  -\frac{1}{2}x_{22} & \vline  & 0 & \vline & 0 & \vline & 0 & \vline & -x_{24} & \vline & 0 & \vline & 0 & \vline & 0 & \vline & 0 & \vline & 0 & \vline & 0 & \vline & 0 & \vline width 1pt \\ \hline
 x_{16} & \vline width 1pt & 0 & \vline  & x_{23} & \vline & 0 & \vline & 0 & \vline & 0 & \vline & 0 & \vline & 0 & \vline & 0 & \vline & 0 & \vline & 0 & \vline & 0 & \vline & 0 & \vline width 1pt \\ \hline
 x_{17} & \vline width 1pt & x_{23} & \vline  & 0 & \vline & x_{24} & \vline & 0 & \vline & 0 & \vline & 0 & \vline & 0 & \vline & 0 & \vline & 0 & \vline & 0 & \vline & 0 & \vline & 0 & \vline width 1pt \\ \hline
 x_{18} & \vline width 1pt & 0     & \vline  & x_{24} & \vline & 0 & \vline & 0 & \vline & 0 & \vline & 0 & \vline & 0 & \vline & 0 & \vline & 0 & \vline & 0 & \vline & 0 & \vline &0 & \vline width 1pt \\ \hline
 x_{19} & \vline width 1pt & 0     & \vline  & 0 & \vline & 0 & \vline & 0 & \vline & 0 & \vline & 0 & \vline & 0 & \vline & 0 & \vline & 0 & \vline & 0 & \vline & 0 & \vline & 0 & \vline width 1pt \\ \hline
 x_{20} & \vline width 1pt & 0     & \vline  & 0 & \vline & 0 & \vline & 0 & \vline & 0 & \vline & 0 & \vline & 0 & \vline & 0 & \vline & 0 & \vline & 0 & \vline & 0 & \vline & 0 & \vline width 1pt \\ \hline
 x_{21} & \vline width 1pt & 0     & \vline  & 0 & \vline & 0 & \vline & 0 & \vline & 0 & \vline & 0 & \vline & 0 & \vline & 0 & \vline & 0 & \vline & 0 & \vline & 0 & \vline & 0 & \vline width 1pt \\ \hline
 x_{22} & \vline width 1pt & 0     & \vline  & 0 & \vline & 0 & \vline & 0 & \vline & 0 & \vline & 0 & \vline & 0 & \vline & 0 & \vline & 0 & \vline & 0 & \vline & 0 & \vline & 0 & \vline width 1pt \\ \hline
 x_{23} & \vline width 1pt & 0     & \vline  & 0 & \vline & 0 & \vline & 0& \vline & 0 & \vline & 0 & \vline & 0 & \vline & 0 & \vline & 0 & \vline & 0 & \vline & 0 & \vline & 0 & \vline width 1pt \\ \hline
 x_{24} & \vline width 1pt & 0     & \vline  & 0 & \vline & 0 & \vline & 0 & \vline & 0 & \vline & 0 & \vline & 0 & \vline & 0 & \vline & 0 & \vline & 0 & \vline & 0 & \vline & 0  & \vline width 1pt\\ \thickhline
\end{array}
$$
\end{landscape}
In the sequel we shall consider $\BB$ as a Lie algebra over an arbitrary field $F$ of characteristic $p$ where $p\neq 2$ ($p=0$ is considerable).
\\ \indent
Consider the symmetric algebra $S(\fn)$ as the polynomial algebra in the variables $x_1,\ldots,x_{24}$ over~$F$. 
We shall expos four invariants $c_1$, $c_2$, $c_3$, $c_4$ and a subset $\{v_{i} \ | \ i\in\{1,\ldots,24\}\setminus\{2,9,16,24\}\}$ of $S(\fn)$ such that 
\begin{equation}
\ad{v_{i}(x_j)}=\left\{\begin{array}{ll}
								0 & i<j \\
								c_{1} & i=j\not\in\{3,4,6,7,10,13\} \\
								c_{2} & i=j\in\{4,7,10,13\} \\
								c_{3} & i=j\in\{3,6\}
							\end{array}\right.
\end{equation}
for $i\leq j$ and  $i,j\not\in\{2,9,16,24\}$.
\\ \indent
Denote $c_{1}=x_{24}$, $c_{2}=2x_{16}x_{24}-2x_{18}x_{23}-x_{20}x_{22}+x_{21}^{2}$.
It is easy to verify that $c_1,c_2\in S(\fn)^{\fn}$, that is, $\ad{x_i}(c_j)=0$ for all $i=1,\ldots,24$, \ $j=1,2$.
\\ \indent
Denote $v_{4}=-2x_{13}x_{24}+2x_{15}x_{23}+x_{17}x_{22}-x_{19}x_{21}$. 
Then $\ad{x_4}(v_4)=-c_2$, \ $\ad{x_i}(v_4)=0$ for $i\neq 4$. 
Denote $u_{9}=x_{9}x_{24}-x_{11}x_{23}+x_{14}x_{22}-\frac{1}{2}x_{19}^{2}$. 
Then $\ad{x_4}(u_9)=v_4$, \ $\ad{x_{i}}(u_9)=0$ for $i\neq 4$. 
It follows that the element $c_{3}=c_{2}u_{9}+\frac{1}{2}v_{4}^{2}$ is invariant. 
That is, $c_3\in S(\fn)^{\fn}$.
\\ \indent
Denote $v_{7}=2x_{10}x_{24}-2x_{12}x_{23}+x_{19}x_{20}-x_{17}x_{21}$. 
Then $\ad{x_3}(v_7)=-v_4$, \ $\ad{x_7}(v_7)=-c_2$, \ $\ad{x_i}(v_7)=0$ for $i\neq 3,7$. 
Denote $u_{6}=x_{6}x_{24}-x_{8}x_{23}+x_{14}x_{21}-\frac{1}{2}x_{17}x_{19}$. 
Then $\ad{x_3}(u_6)=-u_9$, \ $\ad{x_4}(u_6)=-\frac{1}{2}v_7$, \ $\ad{x_7}(u_6)=-\frac{1}{2}v_4$, \ $\ad{x_i}(u_6)=0$ for $i\neq 3,4,7$. 
Denote $v_{3}=-c_{2}u_{6}+\frac{1}{2}v_{4}v_{7}$. 
We have $\ad{x_3}(v_3)=-c_3$, \ $\ad{x_i}(v_3)=0$ for $i\neq 3$. 
Denote $u_{2}=x_{2}x_{24}-x_{5}x_{23}-\frac{1}{2}x_{14}x_{20}+\frac{1}{4}x_{17}^{2}$. 
We have $\ad{x_3}(u_2)=-u_6$, \ $\ad{x_7}(u_2)=-\frac{1}{2}v_{7}$, \ $\ad{x_i}(u_2)=0$ for $i\neq 3,7$. 
Denote $w_{3}=u_9v_7+u_6v_4$. 
Then $\ad{x_4}(w_3)=v_3$, \ $\ad{x_7}(w_3)=-c_3$, \ $\ad{x_i}(w_3)=0$ for $i\neq 4,7$.
It follows that the element $c_{4}=-u_2c_{3}+\frac{1}{2}u_6v_{3}+\frac{1}{4}v_7w_3$ is invariant. 
That is, $c_4\in S(\fn)^{\fn}$.
\\ \indent
We have then four homogeneous invariants $c_1,c_2,c_3,c_4$ and  elements $v_3$, $v_4$, $v_7$.
We shall define the remaining $v_{i}$.
\\ \indent
$v_{23}=x_{1}$, $v_{22}=x_{5}$, $v_{21}=x_{8}$, $v_{20}=-\frac{1}{2}x_{11}$, $v_{19}= x_{12}$, $v_{18}=-x_{14}$, $v_{17}=-x_{15}$, 
$v_{15}=x_{17}$, $v_{14}=x_{18}$, $v_{12}=-x_{19}$, $v_{11}=\frac{1}{2}x_{20}$, $v_{8}=x_{21}$, $v_{5}=-x_{22}$, $v_{1}=-x_{23}$. \\
$v_{13}=2x_{4}x_{24}-2x_{17}x_{18}+2x_{15}x_{20}-2x_{12}x_{21}$. \\ 
$v_{10}=-2x_{7}x_{24}+2x_{18}x_{19}+2x_{12}x_{22}-2x_{15}x_{21}$. 
Note that $\ad{x_2}(v_{10})=-v_7$, \ $\ad{x_6}(v_{10})=-v_4$, \ $\ad{x_{10}}(v_{10})=-c_2$, \ $\ad{x_i}(v_{10})=0$ for $i\neq 2,6,10$. 
\\ \indent  
Let $u_{3}=-x_{3}x_{24}-x_{11}x_{21}+x_{8}x_{22}-x_{15}x_{19}$. 
Then $\ad{x_2}(u_{3})=-u_6$, \ $\ad{x_4}(u_{3})=-\frac{1}{2}v_{10}$, \ $\ad{x_{6}}(u_{3})=u_9$, \ $\ad{x_{10}}(u_{3})=-\frac{1}{2}v_4$, \ $\ad{x_i}(u_{3})=0$ for $i\neq 2,4,6,10$.
Finally $v_{6}=-c_{2}u_{3}+\frac{1}{2}v_{4}v_{10}$. 
\vskip 0.1cm \indent
Assume $p>2$.
We denote by $S_{p}(\fn)$ the polynomial subalgebra of $S(\fn)$ which generated over $F$ by $c_1$ and the $p$-th powers $x_{i}^{p}$: \
$S_{p}(\fn)=F[x_{1}^{p},\ldots,x_{23}^{p},c_{1}]$.
Clearly, $S_{p}(\fn)\subseteq S(\fn)^{\fn}$.
Denote by $Q(A)$ the quotient field of a domain $A$.
The field extension $Q( S_{p}(\fn)[c_2,c_3,c_4])/Q( S_{p}(\fn))$ is of degree $p^3$.
Indeed, $\ad{h_4}(c_1)=0$ and $\ad{h_4}(c_2)=2c_{2}$, thus $c_2\not\in Q(S_{p}(\fn))$.
Since $c_{2}^{p}\in Q(S_{p}(\fn))$ we deduce $[Q( S_{p}(\fn)[c_2]):Q( S_{p}(\fn))]=p$. 
Also, $\ad{h_3}(c_1)=\ad{h_3}(c_2)=0$ and $\ad{h_3}(c_3)=2c_3$ (since $\ad{h_3}(u_9)=2u_9$, $\ad{h_3}(v_4)=v_4$), 
thus $c_3\not\in Q( S_{p}(\fn)[c_2])$. 
Since $c_{3}^{p}\in Q( S_{p}(\fn)[c_2])$ we have $[Q( S_{p}(\fn)[c_2,c_3]):Q( S_{p}(\fn)[c_2])]=p$. 
Now, $\ad{h_2}(c_1)=\ad{h_2}(c_2)=0$, $\ad{h_2}(c_3)=0$ (since $\ad{h_2}(u_9)=\ad{h_2}(v_4)=0$) and $\ad{h_2}(c_4)=2c_4$ (since $\ad{h_2}(u_2)=2u_2$, 
$\ad{h_2}(u_6)=u_6$, $\ad{h_2}(v_3)=v_3$, $\ad{h_2}(v_7)=v_7$, $\ad{h_2}(w_3)=w_3$), thus $c_4\not\in Q( S_{p}(\fn)[c_2,c_3])$. 
Since $c_{4}^{p}\in Q( S_{p}(\fn)[c_2,c_3])$ we deduce $[Q( S_{p}(\fn)[c_2,c_3,c_4]):Q( S_{p}(\fn)[c_2,c_3])]=p$. 
The extension $Q( S_{p}(\fn)[c_2,c_3,c_4])/Q( S_{p}(\fn))$ is therefore of degree $p^3$ as required.
\\ \indent
Using $(1)$, similar considerations yield that the extension $Q(S(\fn)^{\fn}[x_{i} \ | \ i\not\in\{2,9,16,24\}])/Q(S(\fn)^{\fn})$ is of degree $p^{20}$.
The field $Q(S(\fn))=F(x_1,\ldots,x_{24})$ is an extension of degree $p^{23}$ of $Q(S_{p}(\fn))$. 
By degree consideration $Q(S(\fn)^{\fn})=Q(S_{p}(\fn)[c_2,c_3,c_4])$, that is, the domains $S(\fn)^{\fn}$ and  $S_{p}(\fn)[c_2,c_3,c_4]$ have the same quotient field. 
\\ \indent  
The goal is to prove equality of the rings: $S(\fn)^{\fn}=S_{p}(\fn)[c_2,c_3,c_4]$.
The ring $S(\fn)^{\fn}$ is integral over $S_{p}(\fn)[c_2,c_3,c_4]$. 
Therefore, it suffices to prove that $S_{p}(\fn)[c_2,c_3,c_4]$ is normal (integrally closed in its quotient field), equivalently, $S_{p}(\fn)[c_2,c_3,c_4]$ satisfies $(S_1)$ and $(R_2)$ (see [8, p.183]).
Consider the polynomial ring $R=S_{p}(\fn)[t_2,t_3,t_4]$ in the variables $t_2,t_3,t_4$, and its elements $f_i=t_{i}^{p}-c_{i}^{p}$, $i=2,3,4$.
We have $R/(f_2)\cong S_{p}(\fn)[c_2,t_3,t_4]$, $R/(f_2,f_3)\cong S_p(\fn)[c_2,c_3,t_4]$ and $R/(f_2,f_3,f_4)\cong S_p(\fn)[c_2,c_3,c_4]$ (we use [1, Lemma 1.15]) .
Hence, $f_2,f_3,f_4$ form an $R$-sequence. 
$S_p(\fn)[c_2,c_3,c_4]$ is therefore a complete intersection ring (see [2, 1.11].
Moreover, from [5, 2.1.28] it follows that $S_p(\fn)[c_2,c_3,c_4]$ is a Cohen-Macaulay ring, hence satisfies $(S_1)$.
To prove $(R_1)$, we have to show that if $P$ is an element of the singular locus of $R$ such that $f_2,f_3,f_4\in P$, then $\Ht{P}>4$. 
So let $P$ be such a prime.
We have 
$$
\begin{array}{ll}
	\det(\partial(f_2,f_3,f_4)/\partial(x_{16}^{p},x_{9}^{p},x_{2}^{p}))=2c_{1}^{3p}c_{2}^{p}c_{3}^{p}, & 
	\det(\partial(f_2,f_3,f_4)/\partial(x_{16}^{p},x_{9}^{p},x_{6}^{p}))=-2c_{1}^{3p}c_{2}^{p}v_{3}^{p}, \\ \\
	\det(\partial(f_2,f_3,f_4)/\partial(x_{16}^{p},x_{13}^{p},x_{2}^{p}))=-4c_{1}^{3p}v_{4}^{p}c_{3}^{p}, &
     \det(\partial(f_2,f_3,f_4)/\partial(x_{18}^{p},x_{15}^{p},x_{8}^{p}))=-4x_{23}^{3p}v_{4}^{p}v_{3}^{p}.
\end{array}
$$

Therefore, $P$ contains two elements $a,b\in\{c_1^{p},x_{23}^{p},c_2^{p},v_4^{p},c_3^{p},v_3^{p}\}\subseteq S_{p}(\fn)$ hence contains the prime ideal $(a,b,f_2,f_3,f_4)$ which is of height $5$.
We deduce that $S(\fn)^{\fn}=S_{p}(\fn)[c_2,c_3,c_4]$. 
Finally, if $p=0$ then  $S(\fn)^{\fn}=F[c_1,c_2,c_3,c_4]$ by [2, section 3].
\\ \indent
Let $Z(\fn)$ be the center of the enveloping algebra $U(\fn)$.
We shall use the same notation $x_i$, $i=1,\ldots,24$ for the basis of $\fn$, consider it as a Lie subalgebra of $U(\fn)$.  
The elements in $Z(\fn)$ correspond to $c_1$,$c_2$,$c_3$,$c_4$ will respectively denoted by $z_1$,$z_2$,$z_3$,$z_4$. 
Suppose $p>2$.
The analogous polynomial ring to $S_{p}(\fn)$ in $Z(\fn)$ is $Z_{p}(\fn)=F[x_{1}^{p},\ldots,x_{23}^{p},z_{1}]$. 
$S(\fn)$ is isomorphic to the graded algebra of $U(\fn)$, and we set $S(\fn)=\gr{U(\fn)}$.
In particular $c_i=\gr{z_i}$, $x_{i}^{p}=\gr{x_{i}^{p}}=(\gr{x_{i}})^{p}$.
Therefore $S(\fn)^{\fn}=\gr(Z_{p}(\fn)[z_2,z_3,z_4])\subseteq\gr{Z(\fn)}$.
The inclusion $\gr{Z(\fn)}\subseteq S(\fn)^{\fn}$ is trivial. 
Since $S(\fn)^{\fn}=S_{p}(\fn)[c_2,c_3,c_4]$, from [3, p.180, Prop. 10(ii), section 2.9] we have $Z(\fn)=Z_{p}(\fn)[z_2,z_3,z_4]$.
\\ \indent
Let $\varphi:R\to Z(\fn)$ be the $F$-algebra epimorphism defined by $\varphi(x_{i}^{p})=x_{i}^{p}$, \ $\varphi(c_1)=z_1$, \
$\varphi(t_i)=z_{i}$, $i=2,3,4$. 
Obviously, $Rf\subseteq\ker\varphi$. Hence $Z(\fn)$ is a homomorphic image of $R/Rf$. 
The rings $Z(\fn)$ and $R/Rf$ are both domains with equal Krull dimension, hence $Z(\fn)\cong R/Rf$. 
We deduce that $Z(\fn)\cong S(\fn)^{\fn}$. 
The rings $Z(\fn), S(\fn)^{\fn}$ are complete intersection.
\\ \indent
If $p=0$, we have $S(\fn)^{\fn}=F[c_1,c_2,c_3,c_4]=F[\gr{z_1},\gr{z_2},\gr{z_3},\gr{z_4}]=\gr{F[z_1,z_2,z_3,z_4]}\subseteq\gr{Z(\fn)}$ and $\gr{Z(\fn)}\subseteq S(\fn)^{\fn}$.
Therefore $Z(\fn)=F[z_1,z_2,z_3,z_4]$ and $Z(\fn)\cong S(\fn)^{\fn}$ as polynomial algebras in four variables. 
We should remark here that an isomorphism $Z(\fn)\cong S(\fn)^{\fn}$ where $p=0$ in known [6, Proposition 4.8.12]. 
\\ \indent
Let summarize the main results that presented until now: 
\vskip 0.1cm\noindent 
\textbf{2.1. Theorem.} \textbf{a.} \ \textit{Suppose} $\mathit{p>2}$. \textit{Then} 
$$
\mathit{S(\fn)^{\fn}=S_{p}(\fn)[c_2,c_3,c_4]\cong S_{p}(\fn)[t_2,t_3,t_4]/(t_2^{p}-c_{2}^{p},t_3^{p}-c_{3}^{p},t_4^{p}-c_{4}^{p})\cong Z(\fn)=Z_{p}(\fn)[z_2,z_3,z_4]}.
$$
\textit{In particular,} $\mathit{S(\fn)^{\fn}}$, $\mathit{Z(\fn)}$ \textit{are complete intersection rings.} \\ \indent
\textbf{b.} \ \textit{For} $\mathit{p=0}$, \ $\mathit{S(\fn)^{\fn}=F[c_1,c_2,c_3,c_4]\cong Z(\fn)=F[z_1,z_2,z_3,z_4]}$. \\
\textit{In particular,} $\mathit{S(\fn)^{\fn}}$, $\mathit{Z(\fn)}$ \textit{are polynomial rings in four variables.}
\vskip 0.1cm \indent
Let $S(\BB)_{\si}^{\BB}$ be the Poisson semi-center of $S(\BB)$. 
Assume $p>2$.
The linear transformations $\ad{x_i}, \ad{h_{j}}:S(\BB)\to S(\BB)$ satisfy the split equation $X^{p}-X=0$ over $F$. 
Together with the fact $\fn=[\BB,\BB]$, we have $S(\BB)_{\si}^{\BB}=S(\BB)^{\fn}$ (see [2, section 4]). 
Denote $S_{p}(\BB)=S_{p}(\fn)[h_{1}^{p},h_{2}^{p},h_{3}^{p},h_{4}^{p}]$. 
$S_{p}(\BB)$ is a polynomial subalgebra of $S(\BB)$.
Clearly, $S_{p}(\BB)[c_2,c_3,c_4]\subseteq S(\BB)^{\fn}$.
The field extension $Q(S_{p}(\BB)[c_2,c_3,c_4])/Q( S_{p}(\BB))$ is of degree $p^3$. 
The field extension 
$$
	Q(S(\BB)^{\fn}[\{x_{i} \ | \ i\not\in\{2,9,16,24\}\}\cup\{h_1,h_2,h_3,h_4\}])/Q(S(\BB)^{\fn})
$$ 
is of degree $p^{24}$.
Indeed, as before, using $(1)$ we get that the field extension 
$$
	Q(S(\BB)^{\fn}[x_{i} \ | \ i\not\in\{2,9,16,24\}])/Q(S(\BB)^{\fn})
$$ 
is of degree $p^{20}$.
Now, $\ad{c_1}(h_1)=-c_1$ while $\ad{c_1}(x_{i})=0$, $\ad{c_4}(h_2)=-2c_4$ while $\ad{c_4}(x_{i})=\ad{c_4}(h_{1})=0$,
$\ad{c_3}(h_3)=-2c_3$ while $\ad{c_3}(x_{i})=\ad{c_3}(h_{1})=\ad{c_3}(h_{2})=0$, 
$\ad{c_2}(h_4)=-2c_2$ while $\ad{c_2}(x_{i})=\ad{c_2}(h_{1})=\ad{c_2}(h_{2})=\ad{c_2}(h_{3})=0$.
Therefore, the field extension 
$$
	Q(S(\BB)^{\fn}[\{x_{i} \ | \ i\not\in\{2,9,16,24\}\}\cup\{h_1,h_2,h_3,h_4\}])/Q(S(\BB)^{\fn}[x_{i} \ | \ i\not\in\{2,9,16,24\}])
$$
is of degree $p^4$, as required. 
\\ \indent
The field $Q(S(\BB))=F(x_1,\ldots,x_{24},h_1,h_2,h_3,h_4)$ is an extension of degree $p^{27}$ of $Q(S_{p}(\BB))$. 
By degree consideration $Q(S(\BB)^{\fn})=Q(S_{p}(\BB)[c_2,c_3,c_4])$, that is, the domains $S(\BB)^{\fn}$ and  $S_{p}(\BB)[c_2,c_3,c_4]$ have the same quotient field.  
\\ \indent
Identical arguments we applied to $S(\fn)^{\fn}$ yield 
$$
	S(\BB)^{\fn}=S_{p}(\BB)[c_2,c_3,c_4]\cong S_{p}(\BB)[t_2,t_3,t_4]/(t_{2}^{p}-c_{2}^{p},t_{3}^{p}-c_{3}^{p},t_{4}^{p}-c_{4}^{p}),
$$ 
where $t_2,t_3,t_4$ are algebraically independent over $S_{p}(\BB)$.
Also, if $p=0$ then $S(\BB)^{\fn}=F[c_1,c_2,c_3,c_4]$.  
\\ \indent
Let $Sz(\BB)$ be the semi-center of $U(\BB)$. 
By [4, Proposition 2.1], $Sz(\BB)$ is commutative.
We shall use the same notation $h_1$, $h_2$, $h_3$, $h_4$ for the basis of the Cartan subalgebra of $\BB$, consider it as a Lie subalgebra of $U(\BB)$.
Assume $p>2$.
The analogous polynomial ring to $S_{p}(\BB)$ in $Sz(\BB)$ is $Z_{p}(\BB)=Z_{p}(\fn)[h_{1}^{p}-h_{1},h_{2}^{p}-h_{2},h_{3}^{p}-h_{3},h_{4}^{p}-h_{4}]$ ($x_{i}^{p}$ and $h_{j}^{p}-h_{j}$ are central (weight $0$) while $z_1$ is semi-central with non zero weight). 
$S(\BB)$ is isomorphic to the graded algebra of $U(\BB)$, and we set $S(\BB)=\gr{U(\BB)}$.
In particular $h_{i}^{p}=\gr(h_{i}^{p}-h_{i})$.
Therefore $S(\BB)_{\si}^{\BB}=S(\BB)^{\fn}=\gr(Z_{p}(\BB)[z_2,z_3,z_4])\subseteq\gr{Sz(\BB)}$ ($z_2$, $z_3$ and $z_4$ are clearly semi-central with non zero weight).
The inclusion $\gr{Sz(\BB)}\subseteq S(\BB)_{\si}^{\BB}$ is trivial. 
Since $S(\BB)_{\si}^{\BB}=S_{p}(\BB)[c_2,c_3,c_4]$, from [3, p.180, Prop.10(ii), section 2.9] we have $Sz(\BB)=Z_{p}(\BB)[z_2,z_3,z_4]$.
\\ \indent
Let $\varphi:S_{p}(\BB)[t_2,t_3,t_4]\to Sz(\BB)$ be the $F$-algebra epimorphism defined by $\varphi(x_{i}^{p})=x_{i}^{p}$, \ $\varphi(h_{j}^{p})=h_{j}^{p}-h_{j}$ \ $\varphi(c_1)=z_1$, \
$\varphi(t_i)=z_{i}$, $i=2,3,4$. 
Obviously, $t_{i}^{p}-c_{i}^{p}\in\ker\varphi$ for all $i=2,3,4$. 
Hence $Sz(\BB)$ is a homomorphic image of $S_{p}(\BB)[t_2,t_3,t_4]/(t_{2}^{p}-c_{2}^{p},t_{3}^{p}-c_{3}^{p},t_{4}^{p}-c_{4}^{p})$. 
The rings $Sz(\BB)$ and $S_{p}(\BB)[t_2,t_3,t_4]/(t_{2}^{p}-c_{2}^{p},t_{3}^{p}-c_{3}^{p},t_{4}^{p}-c_{4}^{p})$ are both domains with equal Krull dimension, hence $Sz(\BB)\cong S_{p}(\BB)[t_2,t_3,t_4]/(t_{2}^{p}-c_{2}^{p},t_{3}^{p}-c_{3}^{p},t_{4}^{p}-c_{4}^{p})$. 
We deduce that $Sz(\BB)\cong S(\BB)_{\si}^{\BB}$. 
The rings $Sz(\BB), S(\BB)_{\si}^{\BB}$ are complete intersection.
\\ \indent
If $p=0$, we have $S(\BB)_{\si}^{\BB}=F[c_1,c_2,c_3,c_4]=F[\gr{z_1},\gr{z_2},\gr{z_3},\gr{z_4}]=\gr{F[z_1,z_2,z_3,z_4]}\subseteq\gr{Sz(\BB)}$ and $\gr{Sz(\BB)}\subseteq S(\BB)_{\si}^{\BB}$.
Therefore $Sz(\BB)=F[z_1,z_2,z_3,z_4]$ and $Sz(\BB)\cong S(\BB)_{\si}^{\BB}$ as polynomial algebras in four variables. 
We should remark here that an isomorphism $Sz(\BB)\cong S(\BB)_{\si}^{\BB}$ where $p=0$ and $F$ is algebraically closed is known [9]. 
\\ \indent
Let summarize the results for the semi-centers: 
\vskip 0.1cm\noindent 
\textbf{2.2. Theorem.} \textbf{a.} \ \textit{Suppose} $\mathit{p>2}$. \textit{Then}
$$
	\mathit{S(\BB)_{\si}^{\BB}=S_{p}(\BB)[c_2,c_3,c_4]\cong S_{p}(\BB)[t_2,t_3,t_4]/(t_{2}^{p}-c_{2}^{p},t_{3}^{p}-c_{3}^{p},t_{4}^{p}-c_{4}^{p})\cong Sz(\BB)=Z_{p}(\BB)[z_2,z_3,z_4]}.
$$
\textit{The rings} $\mathit{S(\BB)_{\si}^{\BB}}$, $\mathit{Sz(\BB)}$ \textit{are complete intersection.} \\ \indent
\textbf{b.} \ \textit{For} $\mathit{p=0}$, 
$$
	\mathit{S(\BB)_{\si}^{\BB}=F[c_1,c_2,c_3,c_4]\cong Sz(\BB)=F[z_1,z_2,z_3,z_4]}.
$$ 
\textit{Therefore,} $\mathit{S(\BB)_{\si}^{\BB}=S(\fn)^{\fn}}$, \ $\mathit{Sz(\BB)=Z(\fn)}$
\textit{and the both rings} $\mathit{S(\BB)_{\si}^{\BB}}$, $\mathit{Sz(\BB)}$ \textit{are polynomial rings in four variables.}
\vskip 0.1cm \indent
equipped with the semi-centers we can find the Poisson center $S(\BB)^{\BB}$ of $S(\BB)$ and the center $Z(\BB)$ of $U(\BB)$ (reversing the approach in [2]).
\\ \indent
Suppose $p>2$.
Clearly, $Q(F[x_{1}^{p},\ldots,x_{24}^{p},h_{1}^{p},h_{2}^{p},h_{3}^{p},h_{4}^{p}])\subseteq Q(S(\BB)^{\BB})$.
The field extension $Q(S(\BB)^{\BB}[c_1,c_2,c_3,c_4])/Q(S(\BB)^{\BB})$ is of degree $p^{4}$. 
But $S(\BB)^{\BB}[c_1,c_2,c_3,c_4])\subseteq S(\BB)_{\si}^{\BB}=S_{p}(\BB)[c_2,c_3,c_4]$, and $Q(S_{p}(\BB)[c_2,c_3,c_4])$ is of degree $p^{4}$ over $Q(F[x_{1}^{p},\ldots,x_{24}^{p},h_{1}^{p},h_{2}^{p},h_{3}^{p},h_{4}^{p}])$.
By degree consideration we must have $Q(S(\BB)^{\BB})=Q(F[x_{1}^{p},\ldots,x_{24}^{p},h_{1}^{p},h_{2}^{p},h_{3}^{p},h_{4}^{p}])$, hence $S(\BB)^{\BB}=F[x_{1}^{p},\ldots,x_{24}^{p},h_{1}^{p},h_{2}^{p},h_{3}^{p},h_{4}^{p}]$. 
That is, $S(\BB)^{\BB}$ is generated over $F$ by the $p$-th powers of the basis elements of $\BB$.
For $p=0$ we of course have $S(\BB)^{\BB}=F$ by [2, section 3].
\\ \indent
Applying grading consideration one get $Z(\BB)=F[x_{1}^{p},\ldots,x_{24}^{p},h_{1}^{p}-h_1,h_{2}^{p}-h_2,h_{3}^{p}-h_3,h_{4}^{p}-h_4]$ for $p>2$ and $Z(\BB)=F$ for $p=0$. 
\\ \indent
So we have
\vskip 0.1cm\noindent 
\textbf{2.3. Theorem.} \textbf{a.} \ \textit{Suppose} $\mathit{p>2}$. \textit{Then} 
$$
	\mathit{S(\BB)^{\BB}=F[x_{1}^{p},\ldots,x_{24}^{p},h_{1}^{p},h_{2}^{p},h_{3}^{p},h_{4}^{p}]\cong Z(\BB)=F[x_{1}^{p},\ldots,x_{24}^{p},h_{1}^{p}-h_1,h_{2}^{p}-h_2,h_{3}^{p}-h_3,h_{4}^{p}-h_4]},
$$ 
\textit{and are polynomial rings in 28 variables.} \\ \indent
\textbf{b.} \ \textit{For} $\mathit{p=0}$, \ $\mathit{S(\BB)^{\BB}=Z(\BB)=F}$.
\subsection*{3. \ The $\mathbf{C_n}$ case}

\ \ \ \ \ In this case the Borel subalgebra $\BB$ is of dimension $n^2+n$ over $F$, where $F$ is an arbitrary field of characteristic $p\neq 2$  ($p=0$ is considerable).
Its nil radical $\fn$ is of dimension $n^2$ over F consisting of $2n\times 2n$ matrices with its standard basis
$e_{i,j}-e_{n+j,n+i}$ $(1\leq i<j\leq n)$, $e_{i,n+i}$ $(1\leq i\leq n)$ and $e_{i,n+j}+e_{j,n+i}$ $(1\leq i<j\leq n)$.
The standard basis of the Cartan subalgebra is $h_{i}=e_{i,i}-e_{n+i,n+i}$ $(1\leq i\leq n)$. 
\\ \indent
For each positive integer $l$ and for each $s=1,\ldots,l$, let
$$
\fr_{l}(s)=l-s+1.
$$

We arrange the standard basis of $\fn$ in a $2n\times 2n$ matrix $M=[m_{i,j}]$ over the symmetric algebra $S(\fn)$,
symmetrically with respect to its anti-diagonal, in the following way: \ for $1\leq i\leq j\leq n$,
\newpage
\hskip 5.2cm $m_{i,j}=e_{i,j}-e_{n+j,n+i}$ \ $(i\neq j)$. \vskip 0.15cm
\hskip 5.2cm $m_{i,\fr_{2n}(j)}=e_{i,n+j}+e_{j,n+i}$. \vskip 0.15cm
\hskip 5.2cm $m_{j,\fr_{2n}(i)}=m_{i,\fr_{2n}(j)}$.  \vskip 0.15cm
\hskip 5.2cm $m_{\fr_{2n}(j),\fr_{2n}(i)}=m_{i,j}$. \vskip 0.15cm
\hskip 5.2cm Finally, set zeros in the remaining entries.
\vskip 0.1cm
For each $i=1,\ldots,n$, denote by $C_i$ the $i$-th right upper block of the matrix $M$ and by $c_{i}$ its determinant;
\\
$$
c_i=\det(C_i)=\left|\begin{array}{cccc}
            m_{1,\fr_{2n}(i)} & m_{1,\fr_{2n}(i-1)} & \ldots & m_{1,\fr_{2n}(1)} \\
            m_{2,\fr_{2n}(i)} & m_{2,\fr_{2n}(i-1)} & \ldots & m_{2,\fr_{2n}(1)} \\
            \vdots      & \vdots      &        & \vdots  \\
            m_{i,\fr_{2n}(i)} & m_{i,\fr_{2n}(i-1)} & \ldots & m_{i,\fr_{2n}(1)}\end{array}\right|
$$
\vskip 0.1cm
Assume $p>2$.
We denote by $S_{p}(\fn)$ the polynomial subalgebra of $S(\fn)$ which generated over $F$ by $c_1$ and the $p$-th powers of the standard basis of $\fn$.
Let $t_2,\ldots,t_n$ be algebraically independent over $S_{p}(\fn)$.
The Poisson center of $S(\fn)$ is denoted by  $S(\fn)^{\fn}$.
The center of the enveloping algebra $U(\fn)$ will denoted by $Z(\fn)$.
We shall use the same notation for the standard basis of $\fn$, consider it as a Lie subalgebra of $U(\fn)$.  
The elements in $Z(\fn)$ correspond to $c_1\ldots,c_n$ will respectively denoted by $z_1,\ldots,z_n$. 
The analogous polynomial ring to $S_{p}(\fn)$ in $Z(\fn)$ will denoted by $Z_{p}(\fn)$. 
Thus $Z_{p}(\fn)$ is generated over $F$ by $z_1$ and the $p$-th powers of the standard basis of $\fn$. 
From [1, sections 2,5] we have 
\vskip 0.1cm \noindent
\textbf{3.1. Theorem.} \textbf{a.} \ \textit{Suppose} $\mathit{p>2}$. \textit{Then} 
$$
\mathit{S(\fn)^{\fn}=S_{p}(\fn)[c_2,\ldots,c_n]\cong S_{p}(\fn)[t_2,\ldots,t_n]/(t_2^{p}-c_{2}^{p},\ldots,t_n^{p}-c_{n}^{p})\cong Z(\fn)=Z_{p}(\fn)[z_2,\ldots,z_n]}.
$$
\textit{In particular,} $\mathit{S(\fn)^{\fn}}$, $\mathit{Z(\fn)}$ \textit{are complete intersection rings.} \\ \indent
\textbf{b.} \ \textit{For} $\mathit{p=0}$, \ $\mathit{S(\fn)^{\fn}=F[c_1,\ldots,c_n]\cong Z(\fn)=F[z_1,\ldots,z_n]}$. \\
\textit{In particular,} $\mathit{S(\fn)^{\fn}}$, $\mathit{Z(\fn)}$ \textit{are polynomial rings in $n$ variables.}
\vskip 0.1cm \indent
Let $S(\BB)_{\si}^{\BB}$ be the Poisson semi-center of $S(\BB)$. 
Assume $p>2$.
Let $S_{p}(\BB)$ be the polynomial subalgebra of $S(\BB)$ which generated over $S_{p}(\fn)$ by $h_1^{p},\ldots, h_{n}^{p}$. 
Let $t_2,\ldots,t_n$ be algebraically independent over $S_{p}(\BB)$.
The semi-center of $U(\BB)$ will denoted by $Sz(\BB)$. 
We shall use the same notation $h_1,\ldots,h_n$ for the standard basis of the Cartan subalgebra of $\BB$, consider it as a Lie subalgebra of $U(\BB)$.
Let $Z_{p}(\BB)$ be the polynomial subalgebra of $Z(\BB)$ which generated over $Z_{p}(\fn)$ by $h_i^{p}-h_i$, $i=1,\ldots,n$. 
\\ \indent
Identical arguments precede theorems 2.2, 2.3, combining calculations in [1, Proposition 2.11 and Theorem 2.17] yield the following results:
\vskip 0.1cm \noindent
\textbf{3.2. Theorem.} \textbf{a.} \ \textit{Suppose} $\mathit{p>2}$. \textit{Then}
$$
	\mathit{S(\BB)_{\si}^{\BB}=S_{p}(\BB)[c_2,\ldots,c_n]\cong S_{p}(\BB)[t_2,\ldots,t_n]/(t_{2}^{p}-c_{2}^{p},\ldots,t_{n}^{p}-c_{n}^{p})\cong Sz(\BB)=Z_{p}(\BB)[z_2,\ldots,z_n]}.
$$
\textit{The rings} $\mathit{S(\BB)_{\si}^{\BB}}$, $\mathit{Sz(\BB)}$ \textit{are complete intersection.} \\ \indent
\textbf{b.} \ \textit{For} $\mathit{p=0}$, \ $\mathit{S(\BB)_{\si}^{\BB}=F[c_1,\ldots,c_n]\cong Sz(\BB)=F[z_1,\ldots,z_n]}$.  
\textit{Therefore,} $\mathit{S(\BB)_{\si}^{\BB}=S(\fn)^{\fn}}$, \ $\mathit{Sz(\BB)=Z(\fn)}$
\textit{and the both rings} $\mathit{S(\BB)_{\si}^{\BB}}$, $\mathit{Sz(\BB)}$ \textit{are polynomial rings in $n$ variables.}
\vskip 0.1cm \noindent
\textbf{3.3. Theorem.} \textbf{a.} \ \textit{Suppose} $\mathit{p>2}$. \textit{The Poisson center} $\mathit{S(\BB)^{\BB}}$ \textit{coincides with the} $\mathit{p}$\textit{-center, that is,} $\mathit{S(\BB)^{\BB}}$ \textit{is generated over} $\mathit{F}$ \textit{by the} $\mathit{p}$\textit{-th powers of the standard basis elements of} $\mathit{\BB}$.
\textit{The center} $\mathit{Z(\BB)}$ \textit{coincides with the extended} $\mathit{p}$\textit{-center, that is,} $\mathit{Z(\BB)}$ \textit{is generated over} $\mathit{F}$ \textit{by the} $\mathit{p}$\textit{-th powers of the standard basis elements of} $\mathit{\fn}$ \textit{and} $\mathit{h_{i}^{p}-h_{i}},$ $\mathit{i=1,\ldots,n.}$
\textit{In particular,} $\mathit{S(\BB)^{\BB}\cong Z(\BB)}$ \textit{as polynomial rings in} $\mathit{n^{2}+n}$ \textit{variables.}
\\ \indent
\textbf{b.} \ \textit{For} $\mathit{p=0}$, 	$\mathit{S(\BB)^{\BB}=Z(\BB)=F}$. 

\end{document}